\documentclass[12pt,a4paper,reqno]{amsart}
\allowdisplaybreaks
\usepackage{amsmath}
\usepackage{amsfonts}
\usepackage{amssymb,amsthm,amsfonts,amsthm,latexsym,enumerate,url,cases}
\numberwithin{equation}{section}
\usepackage{mathrsfs}
\usepackage{hyperref}
\hypersetup{colorlinks=true,citecolor=blue,linkcolor=blue,urlcolor=blue}
\def\N{{\mathbb N}}

\theoremstyle{plain}
\newtheorem{theorem}{Theorem}

\newtheorem{lemma}{Lemma}
\newtheorem{problem}{Problem}
\newtheorem{corollary}{Corollary}

\theoremstyle{definition}

\usepackage{etoolbox}
\makeatletter
\patchcmd{\@settitle}{\uppercasenonmath\@title}{}{}{}
\patchcmd{\@setauthors}{\MakeUppercase}{}{}{}
\patchcmd{\section}{\scshape}{}{}{}
\makeatother

\begin{document}

\title
[Cross representations of additive complements of $r$-th powers]
{Cross representations of additive complements of $r$-th powers}

\author
[Y. Ding, B. Krause, C. S\'andor, Y. Sun and Z. Zhang] {Yuchen Ding, Ben Krause, Csaba S\'andor, Yu-Chen Sun and Zihan Zhang}

\address{(Yuchen Ding$^{1,2}$) $^1$School of Mathematics,  Yangzhou University, Yangzhou 225002, People's Republic of China}
\address{$^2$HUN-REN Alfr\'ed R\'enyi Institute of Mathematics, Budapest, Pf. 127, H-1364 Hungary}
\email{ycding@yzu.edu.cn}

\address{(Ben Krause) School of Mathematics, University of Bristol, Bristol, BS8 1UG, England}
\email{ben.krause@bristol.ac.uk}

\address{(Csaba S\'andor$^{3,4,5}$) $^3$Department of Stochastics, Institute of Mathematics, Budapest University of Technology and Economics, M\H ugyetem rkp. 3., H-1111,
Budapest, Hungary}
\address{$^4$HUN-REN Alfr\'ed R\'enyi Institute of Mathematics, Re\'altanoda utca 13--15., H-1053 Budapest, Hungary}
\address{$^5$MTA–HUN-REN RI Lend\"ulet ``Momentum'' Arithmetic Combinatorics Research Group, Reáltanoda utca 13--15., H-1053 Budapest, Hungary}
\email{sandor.csaba@ttk.bme.hu}

\address{(Yu-Chen Sun) School of Mathematics, University of Bristol, Bristol, BS8 1UG, England}
\email{yuchensun93@163.com}

\address{(Zihan Zhang) School of Mathematical Sciences and LPMC, Nankai University, Tianjin
300071, People’s Republic of China}
\email{2211056@mail.nankai.edu.cn}

\keywords{additive complement, squares, $r$-th powers, Abel's summation, bipartite graphs, multiplication table problem}
\subjclass[2010]{11B13, 11B75}

\begin{abstract}
Let $\mathbb{N}$ be the set of natural numbers and $\mathcal{S}_r=\big\{1^r, 2^r, 3^r,\cdots\big\}$ the set of $r$-th powers, where $r\ge 2$ is a natural number. Let $\mathcal{W}_r$ be an additive complement of $\mathcal{S}_r$ and
$$
f_r(n)=\#\big\{(w,m^r)\in \mathcal{W}_r\times \mathcal{S}_r: n=w+m^r\big\}.
$$
Motivated by a 1993 conjecture of Cilleruelo, we show that 
$$
\sum_{n\le N}f_r(n)-N\gg_r N^{1-\frac{1}{r}}.
$$
Previously, the bound was only proved for $r=2$.
In the case $r=2$, the lower bound above can be made more explicit as 
$$
\sum_{n\le N}f_2(n)-N\gg N^{3/4-o(1)},
$$
which improves the previous bound $N^{1/2}$ due to Ding, Sun, Wang and Xia.
\end{abstract}
\maketitle

\section{Introduction}
Let $\mathbb{N}$ be the set of natural numbers and $\mathcal{S}=\big\{1^2, 2^2, 3^2,\cdots\big\}$ the set of squares. A subset $\mathcal{W}$ of $\mathbb{N}$ is called an {\it additive complement} of $\mathcal{S}$ if all sufficiently large integers can be written as the sum of an element of $\mathcal{W}$ and an element of $\mathcal{S}$. Erd\H os \cite{Erdos} observed that if $\mathcal{W}$ is an additive complement of $\mathcal{S}$, then one clearly has
$$
\liminf_{N\rightarrow\infty}\mathcal{W}\big(N\big)/\sqrt{N}\ge 1,
$$ 
where $\mathcal{W}(N)=\big|\mathcal{W}\cap[1,N]\big|$.
Erd\H os \cite{Erdos} then asked the natural question whether 
$
\liminf_{N\rightarrow\infty}\mathcal{W}\big(N\big)/\sqrt{N}
$ 
is strictly greater than 1. Erd\H os' problem was answered affirmatively by Moser \cite{Moser}, who obtained the explicit lower bound $1.06$. Moser's bound was later improved several times; see Abbott \cite{Abbott}, Balasubramanian \cite{Balasubramanian}, Balasubramanian and Soundararajan \cite{Balasubramanian-2}, Donagi and Herzog \cite{Donagi}, and Ramana \cite{Ramana,Ramana-2}. The best known lower bound is 
\begin{align}\label{eq-1216-1}
\liminf_{N\rightarrow\infty}\mathcal{W}\big(N\big)/\sqrt{N}\ge \frac{4}{\pi}, 
\end{align}
given independently by Cilleruelo \cite{Cilleruelo}, Habsieger \cite{Habsieger}, Balasubramanian and Ramana \cite{B-Rama}. 

Cilleruelo \cite{Cilleruelo} considered the more general setting obtained by replacing squares with $r$-th powers. Let $r\ge 2$ be a natural number and $\mathcal{S}_r$ the set of $r$-th powers. Let $\mathcal{W}_r$ be an additive complement of $\mathcal{S}_r$, and let
$$
f_r(n)=f_{\mathcal{W}_r,\ \mathcal{S}_r}(n):=\#\big\{(w,m^r)\in \mathcal{W}_r\times \mathcal{S}_r: n=w+m^r\big\}.
$$
Cilleruelo proved that
\begin{align}\label{new-gamma-2}
\liminf_{N\rightarrow\infty}\mathcal{W}_r\big(N\big)/N^{1-1/r}\ge \frac{1}{\Gamma\left(2-\frac{1}{r}\right)\Gamma\left(1+\frac{1}{r}\right)}, 
\end{align}
where $\Gamma$ is Euler's gamma function.
Cilleruelo then conjectured that for any additive complement $\mathcal{W}_r$ of $r$-th powers $\mathcal{S}_r$ we have
$$
\sum_{n\le N}f_r(n)-N\ge rN+o(N).
$$
In the special case $r=2$, the problem was later posed again by Ruzsa \cite{Ruzsa2}. 
\begin{problem}[Ruzsa, 2001]\label{pro:2}
Does the set of squares have an additive complement such that
$\sum_{n\le N}f_2(n)\sim N$?
\end{problem}
Ruzsa's starting point for studying this problem came from his solution \cite{Ruzsa} to Erd\H os' problem \cite{Erdos2} on additive complements of powers of two.

Ben Green observed that for 
$
\mathcal{W}=\{w_1<w_2<\cdots\}
$
with
\begin{align}\label{eq-1216-2}
w_n \sim \frac{\pi^2}{16}n^2, \quad \text{or~equivalently} \quad \mathcal{W}(N)\sim \frac{4}{\pi}\sqrt{N}
\end{align}
as $n\rightarrow\infty$ (resp. $N\rightarrow\infty$)
we have
$$
\lim_{N\rightarrow\infty}\frac{1}{N}\sum_{n\le N}f_2(n)=1.
$$
Ben Green then asked Fang (private communication, see the comments in \cite{Chen-Fang}) whether there exists an additive complement $\mathcal{W}$ of $\mathcal{S}$ such that $\mathcal{W}$ satisfies the growth condition of (\ref{eq-1216-2}). Formally,
\begin{problem}[Ben Green, 2017]\label{pro:1}
Does there exist an additive complement $\mathcal{W}=\{w_n\}_{n=1}^{\infty}$ of squares such that 
$w_n \sim \frac{\pi^2 n^2}{16}$ as $n\rightarrow\infty$?
\end{problem}
Motivated by this, Chen and Fang \cite{Chen-Fang} proved that 
\begin{align}\label{eq-1216-3}
\sum_{n\le N}f_2(n)-N\gg N^{1/4}\log N,
\end{align}
provided that $\mathcal{W}$ is an additive complement of $\mathcal{S}$.
As an application, Chen and Fang showed that
$$
\limsup_{n\rightarrow\infty}\frac{\frac{\pi^2}{16}n^2-w_n}{n^{1/2}\log n}\ge \sqrt{\frac{2}{\pi}}\frac{1}{\log 4}
$$
for any additive complement $\mathcal{W}=\{w_n\}_{n=1}^{\infty}$ of $\mathcal{S}$, which was later improved by the first named author \cite{Ding} to
\begin{align*}
\limsup_{n\rightarrow\infty}\frac{\frac{\pi^2}{16}n^2-w_n}{n}\ge \frac{\pi}{4}.
\end{align*}
In a more recent article, Ding, Sun, Wang and Xia \cite{DSWX} improved the bound (\ref{eq-1216-3}) to
\begin{align}\label{eq-1216-5}
\sum_{n\le N}f_2(n)-N\ge 0.193N^{1/2}.
\end{align}
The bounds (\ref{eq-1216-3}) and (\ref{eq-1216-5}) represent progress towards Cilleruelo's conjecture for $r=2$. Unfortunately, for $r\ge 3$, the methods in \cite{Chen-Fang} and \cite{DSWX} that lead to (\ref{eq-1216-3}) and (\ref{eq-1216-5}) do not seem to yield even
$$
\Big(\sum_{n\le N}f_r(n)-N \Big)\rightarrow\infty \quad \text{as }N\rightarrow\infty.
$$
The following new theorem not only extends the bound (\ref{eq-1216-5}) to any $r\ge 2$, but also applies to an even more general setting.

\begin{theorem}\label{thm:1}
Let $r\ge 2$ be a natural number and $\mathscr{S}\subset \mathbb{N}$ be a subset such that
$$
\mathscr{S}(x)\sim c_r x^{1/r} \quad \text{as }x\rightarrow\infty, 
$$
where $c_r>0$ is a constant. Then, for any additive complement $\mathscr{W}$ of $\mathscr{S}$, we have
$$
\sum_{n\le N}f_{\mathscr{S},\mathscr{W}}(n)-N\gg N^{1-\frac{1}{r}},
$$
where the implied constant depends at most on $r$ and $c_r$, and 
$$
f_{\mathscr{S},\mathscr{W}}(n)=\#\left\{(s,w)\in \mathscr{S}\times \mathscr{W}: n=s+w\right\}.
$$
\end{theorem}

The following corollary is an immediate consequence of Theorem \ref{thm:1}.

\begin{corollary}
   Let $\mathcal{W}_r$ be an additive complement of $r$-th powers $\mathcal{S}_r$. Then we have
$$
\sum_{n\le N}f_r(n)-N\gg N^{1-\frac{1}{r}},
$$
where the implied constant depends only on $r$.
\end{corollary}

In \cite{DSWX}, the authors mentioned the goal of obtaining a lower bound of the form 
$$
\sum_{n\le N}f_2(n)-N\gg \rho(N)N^{1/2}
$$ 
with $\rho(N)\rightarrow\infty$ as $N\rightarrow\infty$. Our next theorem not only achieves this goal but also improves the exponent of $N$. 
For an integer \(h\ge1\), let \(\tau(h)\) denote the number of its positive divisors, and put
\[
 T(N)=\max_{1\le h\le N}\tau(h).
\]

For additive complements of the squares, we use $f(n)$ as an abbreviation for $f_2(n)$. 

\begin{theorem}\label{thm:main}
Let \(\mathcal{W}\) be an additive complement of the squares.  Then, for all sufficiently large \(N\),
\[
 \sum_{n\le N}f(n)-N\gg \frac{N^{3/4}}{\sqrt{T(N)}},
\]
where the implied constant depends only on \(\mathcal{W}\).
\end{theorem}

The classical maximal order estimate for the divisor function then gives the following more familiar form.

\begin{corollary}\label{cor:wigert}
For every additive complement \(\mathcal{W}\) of the squares,
\[
 \sum_{n\le N}f_2(n)-N
 \ge N^{3/4}
 \exp\!\left(-\left(\frac{\log 2}{2}+o(1)\right)
       \frac{\log N}{\log\log N}\right).
\]
In particular, for every fixed \(\varepsilon>0\),
\[
 \sum_{n\le N}f_2(n)-N\gg N^{3/4-\varepsilon}.
\]
\end{corollary}
The improved lower bound for $r=2$ comes from the arithmetic features of squares. Thus, the argument of Theorem \ref{thm:main} cannot be applied to the general case $r\ge 3$.

On the Erd\H os Problems website \cite[Problem $\#$33]{Bloom}, Erd\H os is recorded as asking for the smallest possible value of
$$
\limsup_{N\rightarrow\infty}\mathcal{W}\big(N\big)/\sqrt{N}
$$
over all additive complements $\mathcal{W}$ of the squares; denote this value by $L_{\mathcal{S}}$. The website also notes that van Doorn constructed an additive complement $\mathcal{W}_0$ of $\mathcal{S}$ such that
$$
\limsup_{N\rightarrow\infty}\mathcal{W}_0\big(N\big)/\sqrt{N}\le 2\left(\frac{1+\sqrt{5}}{2}\right)^{5/2}\approx 6.66.
$$

\section{Proof of Theorem \ref{thm:1}}
\begin{proof}[Proof of Theorem \ref{thm:1}]
Throughout the proof, let $N$ be sufficiently large.
Let $K$ be a fixed large integer to be chosen later. For any $1\le k\le K$, let
$$
\mathscr{W}_{k,N}=\mathscr{W}\cap \left[\frac{k-1}{2K}N,\ \frac{k}{2K}N\right), \quad h(k)=\big|\mathscr{W}_{k,N}\big|
$$
and
$$
\mathscr{S}_{k,N}=\mathscr{S}\cap \left[\frac{k-1}{2K}N,\ \frac{k}{2K}N\right), \quad g(k)=\big|\mathscr{S}_{k,N}\big|.
$$
Then
\begin{align}\label{eq-1214-2}
g(k)=\Big(c_r+o(1)\Big)\left(\frac{N}{2K}\right)^{1/r}\Big(k^{1/r}-(k-1)^{1/r}\Big).
\end{align}
For any $1\le k\le K$, let
$$
A(k)=\sum_{i=1}^{k}h(i).
$$
Since $\mathscr{W}$ is an additive complement of $\mathscr{S}$, we have
$$
\mathscr{W}(x)\mathscr{S}(x)\ge x+O(1),
$$
which implies that
\begin{align}\label{eq-0320-1}
\mathscr{W}(x)\ge \Big(c_r^{-1}+o(1)\Big)x^{1-1/r},
\end{align}
provided that $x$ is sufficiently large.
Then, by the definition of $h(i)$ and (\ref{eq-0320-1}), we have 
\begin{align}\label{eq-1228-1}
A(k)=\mathscr{W}\left(\frac{kN}{2K}\right)\ge \frac{c_r^{-1}}{1.9}\left(\frac{kN}{2K}\right)^{1-1/r}.
\end{align}

It is easy to see that for any $1\le k\le K$ and $w\in \mathscr{W}_{k,N},\ s\in \mathscr{S}_{k,N}$ we have
$$
w-s\in \left(-\frac{N}{2K}, \frac{N}{2K}\right),
$$
and the number of such differences $w-s$, counted with multiplicity, is
\begin{align}\label{eq-1228-2}
\sum_{1\le k\le K}h(k)g(k)&=\sum_{1\le k\le K}\Big(A(k)-A(k-1)\Big)g(k)\nonumber\\
&=A(K)g(K)+\sum_{1\le k\le K-1}A(k)\Big(g(k)-g(k+1)\Big),
\end{align}
where the last equality follows from Abel summation. Note that
$$
\Big(k^{1/r}-(k-1)^{1/r}\Big)-\Big((k+1)^{1/r}-k^{1/r}\Big)=2k^{1/r}-\Big((k-1)^{1/r}+(k+1)^{1/r}\Big)>0
$$
from Jensen's inequality. In view of (\ref{eq-1214-2}), we know that $g(k)>g(k+1)$ for any $1\le k\le K-1$. Hence, using again the Abel summation formula we deduce from (\ref{eq-1228-1}), (\ref{eq-1228-2}) and (\ref{eq-1214-2}) that
\begin{align*}
\sum_{1\le k\le K}&h(k)g(k)\nonumber\\
&\ge \frac{c_r^{-1}}{1.9}\left(\frac{KN}{2K}\right)^{1-\frac{1}{r}}g(K)+\sum_{1\le k\le K-1}\frac{c_r^{-1}}{1.9}\left(\frac{kN}{2K}\right)^{1-\frac{1}{r}}\Big(g(k)-g(k+1)\Big)\nonumber\\
&=\frac{c_r^{-1}}{1.9}\sum_{1\le k\le K}\left(\left(\frac{kN}{2K}\right)^{1-1/r}-\left(\frac{(k-1)N}{2K}\right)^{1-1/r}\right)g(k)\nonumber\\
&\ge \frac{N}{4K}\sum_{1\le k\le K}\left(k^{1-1/r}-\left(k-1\right)^{1-1/r}\right)\Big(k^{1/r}-(k-1)^{1/r}\Big).
\end{align*}
By the mean value theorem, there exist two numbers $\theta_1$ and $\theta_2$ with $k-1\le \theta_1, \theta_2\le k$ such that
$$
k^{1-1/r}-\left(k-1\right)^{1-1/r}=\left(1-\frac{1}{r}\right)\theta_1^{-1/r}\ge \left(1-\frac{1}{r}\right)k^{-1/r}
$$
and
$$
k^{1/r}-(k-1)^{1/r}=\frac{1}{r}\theta_2^{1/r-1}\ge \frac{1}{r}k^{1/r-1},
$$
from which we conclude
\begin{align}\label{equation-2801}
\sum_{1\le k\le K}h(k)g(k)\ge \frac{N}{4K}\left(1-\frac{1}{r}\right)\frac{1}{r}\sum_{1\le k\le K}\frac{1}{k}\ge \frac{N\log K}{8K}\left(1-\frac{1}{r}\right)\frac{1}{r}.
\end{align}

Let $d\in \left(-\frac{N}{2K}, \frac{N}{2K}\right)$ be an integer attained as a difference $w-s$. Let $\ell_d$ denote the number of pairs $\big(w,s\big)$ such that $d=w-s$. Then by (\ref{equation-2801}) we have
\begin{align*}
\sum_{\substack{-\frac{N}{2K}<d<\frac{N}{2K}\\ \ell_d\ge 1}}\ell_d\ge \frac{N\log K}{8K}\left(1-\frac{1}{r}\right)\frac{1}{r},
\end{align*}
from which it follows that
\begin{align}\label{eq-1214-5}
\sum_{\substack{-\frac{N}{2K}<d<\frac{N}{2K}\\ \ell_d\ge 2}}\ell_d\ge\frac{N\log K}{8K}\left(1-\frac{1}{r}\right)\frac{1}{r}-\frac{N}{K}=\frac{N}{K}\left(\frac{(r-1)\log K}{8r^2}-1\right).
\end{align}

Now, for any $d\in \left(-\frac{N}{2K}, \frac{N}{2K}\right)$ with $\ell_d\ge 2$, we assume that
$$
d=w_1-s_1=w_2-s_2=\cdots=w_{\ell_d}-s_{\ell_d}.
$$
For any $1\le i\neq j\le \ell_d$, we see that
$$
w_i-s_i=w_j-s_j.
$$
It then follows that
\begin{align}\label{eq-1214-6}
w_i+s_j=w_j+s_i<\frac{N}{2}+\frac{N}{2}=N.
\end{align}
For $1\le n\le N$, we define the set $B_n$ to be
\begin{align*}
    B_n=\Big\{ \big\{(w,s), (w',s')  &\big\} :~  n=w+s=w'+s', \\
    &w, w'\in \mathscr{W},\ s,s'\in \mathscr{S},\ w, w', s, s' \le N/2,\ w\neq w' \Big\}.
\end{align*}
Clearly, we have
\[
 |B_n|\leq\binom{f_{\mathscr{S},\mathscr{W}}(n)}{2}.
\]
For any $d\in \left(-\frac{N}{2K}, \frac{N}{2K}\right)$ with $\ell_d\ge 2$ and any $1\le i<j\le \ell_d$, equation (\ref{eq-1214-6}) gives an integer $n\in[1,N]$ such that
\begin{align}\label{eq-0104-1}
    \Big\{(w_i,s_j),(w_j,s_i)\Big\}\in B_n.
\end{align}
Moreover, these elements of the sets $B_n$ are distinct as $d$ and $\{i,j\}$ vary. Indeed, from
$\big\{(w_i,s_j),(w_j,s_i)\big\}$ one recovers the two original pairs
$(w_i,s_i)$ and $(w_j,s_j)$, and hence their common difference $d$.
Therefore,
\begin{align}\label{eq-0104-2}
\sum_{\substack{1\leq n \leq N \\ f_{\mathscr{S}, \mathscr{W}}(n)\ge 2}}\binom{f_{\mathscr{S}, \mathscr{W}}(n)}{2}\geq \sum_{\substack{-\frac{N}{2K}<d<\frac{N}{2K}\\ \ell_d\ge 2}}\binom{\ell_d}{2}.
\end{align}

Now, we deduce from (\ref{eq-0104-2}) that
\begin{align*}
\sum_{\substack{1\leq n \leq N \\ f_{\mathscr{S}, \mathscr{W}}(n)\ge 2}}f_{\mathscr{S}, \mathscr{W}}(n)\big(f_{\mathscr{S}, \mathscr{W}}(n)-1\big)\ge 2\sum_{\substack{-\frac{N}{2K}<d<\frac{N}{2K}\\ \ell_d\ge 2}}\binom{\ell_d}{2},
\end{align*}
which implies
\begin{align}\label{eq-0320-2}
\sum_{\substack{1\leq n \leq N}}f_{\mathscr{S}, \mathscr{W}}(n)\big(f_{\mathscr{S}, \mathscr{W}}(n)-1\big)\ge \sum_{\substack{-\frac{N}{2K}<d<\frac{N}{2K}\\ \ell_d\ge 2}}\ell_d\big(\ell_d-1\big)\ge \sum_{\substack{-\frac{N}{2K}<d<\frac{N}{2K}\\ \ell_d\ge 2}}\ell_d.
\end{align}
Inserting (\ref{eq-1214-5}) into (\ref{eq-0320-2}) gives
\begin{align}\label{eq-0320-3}
\sum_{\substack{1\leq n \leq N}}f_{\mathscr{S}, \mathscr{W}}(n)-N\ge \frac{N}{K}\left(\frac{(r-1)\log K}{8r^2}-1\right) \left(\max_{1\le n\le N}f_{\mathscr{S}, \mathscr{W}}(n)\right)^{-1}.
\end{align}
Trivially, 
$$
\max_{1\le n\le N}f_{\mathscr{S}, \mathscr{W}}(n)\le \mathscr{S}(N)=\big(c_r+o(1)\big)N^{1/r}.
$$
Choosing $K=\lfloor e^{32r}\rfloor+1$, we obtain 
$$
\sum_{\substack{1\leq n \leq N}}f_{\mathscr{S}, \mathscr{W}}(n)-N\gg_{r, c_r} N^{1-1/r}
$$
from (\ref{eq-0320-3}), proving our theorem.
\end{proof}

\section{Proofs of Theorem \ref{thm:main} and Corollary \ref{cor:wigert}}
The proof of Theorem \ref{thm:main} has three ingredients. First, we prove a second-moment lower bound for a restricted representation function. Second, we use a remarkable result on the multiplication table problem. Third, we establish a new lemma on codegrees in bipartite graphs.

\begin{lemma}\label{Lemma-restricted}
Let $\mathcal{S}:=\mathcal{S}_2$ be the set of squares and $\mathcal{W}$ an additive complement of $\mathcal{S}$. Let $N$ be sufficiently large and $K_0=\lfloor e^{128}\rfloor+1$. Put
$$
f^*(n)=\#\left\{(w, m^2)\in \mathcal{W}\times\mathcal{S}: n=w+m^2, m^2\ge \frac{N}{2K_0}\right\}.
$$
Then at least one of the following two estimates holds:
\begin{align}\label{first-moment}
\sum_{n\le N}f(n)-N\ge \frac{\sqrt{2}-1}{2K_0}N,
\end{align}
\begin{align*}
\sum_{\substack{N/K_0\leq n \leq N}}\binom{f^*(n)}{2}\ge \frac{N\log K_0}{128K_0}.
\end{align*}
\end{lemma}
\begin{proof}
Let $N$ be sufficiently large.
Taking $\mathscr{S}=\mathcal{S}$ and $\mathscr{W}=\mathcal{W}$ in Theorem \ref{thm:1}, we have $r=2$ and $c_r=1$.
By (\ref{equation-2801}) with $K=K_0$ we get
\begin{align}\label{eq-0322-1}
\sum_{1\le k\le K_0}h(k)g(k)\ge \frac{N\log K_0}{32K_0},
\end{align}
where
$$
h(k)=\bigg|\mathcal{W}\cap\left[\frac{k-1}{2K_0}N,\ \frac{k}{2K_0}N\right)\bigg| \quad \text{and} \quad g(k)=\bigg|\mathcal{S}\cap\left[\frac{k-1}{2K_0}N,\ \frac{k}{2K_0}N\right)\bigg|.
$$
We split the proof into two cases.

{\it Case I.} 
$$
h(1)\ge 2\sqrt{\frac{N}{K_0}}.
$$
In this case, we clearly have
\begin{align*}
    \sum_{n\le N}f(n)-N\ge \sum_{n\le N/K_0}f(n)-\frac{N}{K_0}\ge h(1)g(1)-\frac{N}{K_0}.
\end{align*}
Since 
$$
g(1)=\big(1+o(1)\big)\sqrt{\frac{N}{2K_0}},
$$
it follows that
\begin{align*}
    \sum_{n\le N}f(n)-N
    \ge \frac{\sqrt{2}-1}{2K_0}N.
\end{align*}

{\it Case II.} 
$$
h(1)< 2\sqrt{\frac{N}{K_0}}.
$$
By (\ref{eq-0322-1}) we obtain
\begin{align*}
\sum_{2\le k\le K_0}h(k)g(k)=\sum_{1\le k\le K_0}h(k)g(k)-h(1)g(1)\ge \frac{N\log K_0}{64K_0}.
\end{align*}
For this restricted part, let \(\ell_d\) denote the number of pairs \((w,s)\) such that, for some \(2\le k\le K_0\),
\[
w\in\mathcal{W}\cap\left[\frac{k-1}{2K_0}N,\frac{k}{2K_0}N\right),
\quad
s\in\mathcal{S}\cap\left[\frac{k-1}{2K_0}N,\frac{k}{2K_0}N\right),
\quad d=w-s.
\]
Then \(d\in(-N/(2K_0),N/(2K_0))\), and the preceding estimate gives
\begin{align*}
\sum_{\substack{-\frac{N}{2K_0}<d<\frac{N}{2K_0}\\ \ell_d\ge 1}}\ell_d
=\sum_{2\le k\le K_0}h(k)g(k)
\ge \frac{N\log K_0}{64K_0}.
\end{align*}
Suppose that \(\ell_d\ge2\), and write the pairs counted by \(\ell_d\) as
\[
d=w_1-s_1=w_2-s_2=\cdots=w_{\ell_d}-s_{\ell_d}.
\]
For \(1\le i<j\le \ell_d\), we have
\[
n=w_i+s_j=w_j+s_i.
\]
Since all the pairs counted here come from blocks with \(k\ge2\), both square parts \(s_i,s_j\) are at least \(N/(2K_0)\). Moreover, all four terms \(w_i,w_j,s_i,s_j\) are smaller than \(N/2\), while \(w_i,s_j\ge N/(2K_0)\). Hence
\[
\frac{N}{K_0}\le n<N,
\]
and the two distinct representations \((w_i,s_j)\) and \((w_j,s_i)\) are counted by \(f^*(n)\). As in the proof of Theorem \ref{thm:1}, these unordered pairs of representations are distinct as \(d\) and \(\{i,j\}\) vary. Consequently,
\begin{align*}
\sum_{\substack{N/K_0\leq n \leq N}}\binom{f^*(n)}{2}
\geq \sum_{\substack{-\frac{N}{2K_0}<d<\frac{N}{2K_0}\\ \ell_d\ge 2}}\binom{\ell_d}{2}.
\end{align*}
Finally, there are at most \(N/K_0+1\) possible integer values of \(d\). Since \(\binom{u}{2}\ge u-1\) for every integer \(u\ge1\), the preceding lower bound for \(\sum \ell_d\) gives
\begin{align*}
\sum_{\substack{-\frac{N}{2K_0}<d<\frac{N}{2K_0}\\ \ell_d\ge 2}}\binom{\ell_d}{2}
&=\sum_{\substack{-\frac{N}{2K_0}<d<\frac{N}{2K_0}\\ \ell_d\ge 1}}\binom{\ell_d}{2}\\
&\ge \sum_{\substack{-\frac{N}{2K_0}<d<\frac{N}{2K_0}\\ \ell_d\ge 1}}(\ell_d-1)\\
&\ge \frac{N\log K_0}{128K_0},
\end{align*}
where the last inequality uses \(\log K_0>128\) and \(N\) sufficiently large. Combining the last two displayed estimates completes the proof of our lemma.
\end{proof}

We will use the upper bound in the following result \cite[Corollary 3]{Ford} on the multiplication table problem.
\begin{lemma}\label{Lemma-table}
Let $L(x)$ be the number of positive integers $n\le x$ which can
be written as $n=m_1m_2$, where each $m_i\le \sqrt{x}$. Then we have
$$
C_1\frac{x}{(\log x)^\delta(\log\log x)^{3/2}}<L(x)<C_2 \frac{x}{(\log x)^\delta(\log\log x)^{3/2}} ,
$$
where $C_1$ and $C_2$ are absolute constants and $\delta=1-\frac{1+\log\log 2}{\log 2}=0.086071\ldots$.
\end{lemma}
For further work on the multiplication table problem, see Erd\H os \cite{Erdos-New} and Tenenbaum \cite{Tenenbaum}. For extensions, see Chang \cite{Chang}, Elekes and Ruzsa \cite{Elekes-Ruzsa}, and Xu and Zhou \cite{Xu-Zhou}.

A \emph{codegree} is the number of common neighbours of two vertices on the same side of a bipartite graph. 
We will need the following graph-theoretic lemma for bipartite graphs with bounded codegrees. Let $\Gamma(x)$ denote the neighbourhood of a vertex $x$.
\begin{lemma}\label{lem:graph}
Let \(G=(X,Y,E)\) be a finite simple bipartite graph.  
Assume that \(|Y|\ge1\), \(T\ge1\), and
\[
 |\Gamma(x)|\le B\quad(x\in X),
\]
and that
\[
 |\Gamma(x)\cap\Gamma(x')|\le T
 \quad\text{whenever }x\ne x'.
\]
Write
\[
 a=\sum_{x\in X}\big(|\Gamma(x)|-1\big)_+,
 \qquad
 Q=\sum_{x\in X}\binom{|\Gamma(x)|}2,
\]
where \((u)_+=\max(u,0)\).  Then
\begin{equation}\label{eq:graph-lemma}
 Q\le a\sqrt{|Y|T}+\frac23B|Y|.
\end{equation}
\end{lemma}
\begin{proof}
We split the vertices of \(X\) according as \(|\Gamma(x)|\le 2\sqrt{|Y|T}\) or \(|\Gamma(x)|>2\sqrt{|Y|T}\).

For a vertex $x$ with \(|\Gamma(x)|\le 2\sqrt{|Y|T}\),
\[
 \binom{|\Gamma(x)|}2
 =\frac{|\Gamma(x)|(|\Gamma(x)|-1)}2
 \le\sqrt{|Y|T}(|\Gamma(x)|-1)_+.
\]
Therefore, the total contribution of these low-degree vertices is at most
\begin{equation}\label{eq:low-degree}
 \sum_{x\in X}\sqrt{|Y|T}(|\Gamma(x)|-1)_+=a\sqrt{|Y|T}.
 \end{equation}

Now we focus on the high-degree vertices of the graph. Set
\[
 X_{\mathrm{h}}=\big\{x\in X:|\Gamma(x)|>2\sqrt{|Y|T}\big\}.
\]
Let
\[
 I=\sum_{x\in X_{\mathrm{h}}}|\Gamma(x)|,
\]
so \(I\) is the number of edges incident with \(X_{\mathrm{h}}\).  For any \(y\in Y\), let \(e_y\) be its degree into \(X_{\mathrm{h}}\).  Then we have
\[
 \sum_{y\in Y}e_y=I.
\]
The Cauchy--Schwarz inequality gives
$$
I^2=\Big(\sum_{y\in Y}e_y\Big)^2\le |Y|\sum_{y\in Y}e_y^2,
$$
or equivalently,
\[
 \sum_{y\in Y}e_y^2\ge\frac{I^2}{|Y|}.
\]
Hence,
\begin{equation}\label{eq:right-pairs-lower}
 \sum_{y\in Y}\binom{e_y}{2}
 =\frac12\left(\sum_{y\in Y}e_y^2-I\right)
 \ge\frac12\left(\frac{I^2}{|Y|}-I\right).
\end{equation}
The left-hand side of \eqref{eq:right-pairs-lower} counts triples \((\{x,x'\},y)\) in which \(x,x'\) are distinct high-degree vertices and \(y\) is a common neighbour.  Since every pair \(x\ne x'\) has codegree at most \(T\) by hypothesis of the lemma,
\begin{equation}\label{eq:right-pairs-upper}
 \sum_{y\in Y}\binom{e_y}{2}
 \le T\binom {|X_{\mathrm{h}}|}2
 \le\frac{T|X_{\mathrm{h}}|^2}{2}.
\end{equation}
Every high-degree vertex has degree greater than \(2\sqrt{|Y|T}\) by our assumption, so
\[
 |X_{\mathrm{h}}|<\frac{I}{2\sqrt{|Y|T}}.
\]
Combining \eqref{eq:right-pairs-lower} and \eqref{eq:right-pairs-upper} yields
\begin{align*}
 \frac{I^2}{|Y|}-I
 &\le T|X_{\mathrm{h}}|^2
 <\frac{TI^2}{4|Y|T}
 =\frac{I^2}{4|Y|}.
\end{align*}
Thus,
\[
 \frac{3I^2}{4|Y|}<I.
\]
If \(I=0\), there is nothing to prove.  Otherwise division by \(I\) gives
\begin{equation*}
 I<\frac{4|Y|}{3}.
\end{equation*}
Finally, because every left degree is at most \(B\),
\begin{align}
 \sum_{x\in X_{\mathrm{h}}}\binom{|\Gamma(x)|}2
 &\le\frac B2\sum_{x\in X_{\mathrm{h}}}|\Gamma(x)|
 =\frac B2I
 <\frac23B|Y|.
 \label{eq:high-degree}
\end{align}
Adding \eqref{eq:low-degree} and \eqref{eq:high-degree} proves \eqref{eq:graph-lemma}.
\end{proof}

Now, we are ready to provide the proof of Theorem \ref{thm:main}.

\begin{proof}[Proof of Theorem \ref{thm:main}]
Let $N$ be sufficiently large. If (\ref{first-moment}) in Lemma \ref{Lemma-restricted} holds, then the theorem follows immediately. From now on, we will assume 
\begin{align}\label{second-moment}
\sum_{\substack{N/K_0\leq n \leq N }}\binom{f^*(n)}{2}\ge \frac{N\log K_0}{128K_0}.
\end{align}
We split the proof into two cases.

{\it Case I.} 
$$
\max_{N/K_0\le n\le N}f^*(n)\ge \frac{25C_2K_0\sqrt{N}}{(\log N)^{\delta}(\log\log N)^{3/2}},
$$
where $\delta$ and $C_2$ are as in Lemma \ref{Lemma-table}.
In this case, there exist 
$$
N/K_0\le n^*\le N
$$ 
and 
$$
1\le m_1<m_2<\cdots<m_v\le \sqrt{N}
$$ 
with 
$$
v\ge \frac{25C_2K_0\sqrt{N}}{(\log N)^{\delta}(\log\log N)^{3/2}}
$$
such that
$$
w_i=n^*-m_i^2\in \mathcal{W},\ 1\le i\le v.
$$
Moreover, $m_i^2\ge \frac{N}{2K_0}$ for any $1\le i\le v$ by the definition of $f^*$.
For any $1\le i\le v$ and $1\le m< m_i$ we have
\begin{align}\label{eq-0320-5}
w_i+m^2=w_i+m_i^2-(m_i^2-m^2)=n^*-(m_i-m)(m_i+m).
\end{align}
Let 
\begin{align*}
U=\{(m_i-m)(m_i+m): 1\le i\le v, 1\le m< m_i\}:=\{u_1<u_2<\cdots<u_t\}.
\end{align*}
Note that
$1\le m_i-m,\ m_i+m \le 2\sqrt{N}$ and 
$
u_t< N.
$
Applying Lemma \ref{Lemma-table} with $x=4N$, we have
\begin{align}\label{eq-0320-6}
t<\frac{5C_2N}{(\log N)^\delta(\log\log N)^{3/2}}.
\end{align}
In view of (\ref{eq-0320-5}) and the definitions of $u_j$ we get
\begin{align}\label{eq-0320-7}
 \sum_{j=1}^{t}f(n^*-u_j)\ge \sum_{i=1}^{v}(m_i-1)\ge \frac{\sqrt{N}}{2\sqrt{K_0}}\sum_{i=1}^{v}1 \ge \frac{10C_2 N}{(\log N)^{\delta}(\log\log N)^{3/2}}.
\end{align}
Thus, by (\ref{eq-0320-6}) and (\ref{eq-0320-7}) we conclude that
\begin{align*}
  \sum_{n\le N}f(n)-N\ge   \sum_{j=1}^{t}f(n^*-u_j)-t>\frac{5C_2N}{(\log N)^\delta(\log\log N)^{3/2}},
\end{align*}
which is stronger than the bound claimed in the theorem.

{\it Case II.} 
\begin{align}\label{degree-1}
\max_{N/K_0\le n\le N}f^*(n)\le \frac{25C_2K_0\sqrt{N}}{(\log N)^{\delta}(\log\log N)^{3/2}}.
\end{align}
In this case, we construct a bipartite graph \(G_N=(X_N,Y_N,E_N)\) as follows. Put
\begin{align*}
 X_N&=\big\{n\in\N:N/K_0<n\le N\big\},\\
 Y_N&=\mathcal{W}\cap\left[1,\Big(1-\frac{1}{2K_0}\Big)N\right],
\end{align*}
and join \(n\in X_N\) to \(w\in Y_N\) if
\[
 n-w=m^2\ge\frac{N}{2K_0}
\]
for some \(m\in\N\).  For fixed \(n,w\), the integer \(m\) is unique, so the graph is simple.  Conversely, every representation counted by \(f^*(n)\) has
\[
 w=n-m^2\le N-\frac{N}{2K_0}=\left(1-\frac{1}{2K_0}\right)N.
\]
Therefore the left degree $|\Gamma(n)|$ of each $n$ is exactly $f^*(n)$.
By \eqref{degree-1}, every left degree $|\Gamma(n)|$ satisfies
\begin{equation}\label{eq:graph-B}
 |\Gamma(n)|\le\frac{25C_2K_0\sqrt{N}}{(\log N)^{\delta}(\log\log N)^{3/2}}.
\end{equation}

We next bound the size
\(
 |Y_N|.
\)
For any \(w\in Y_N\) and any integer \(1\le m\le\left\lfloor\sqrt{\frac{N}{2K_0}}\right\rfloor\),
\[
 w+m^2\le\left(1-\frac{1}{2K_0}\right)N+\frac{N}{2K_0}=N.
\]
All these pairs are counted by \(\sum_{n\le N}f(n)\).  Hence,
\[
 |Y_N|\left\lfloor\sqrt{\frac{N}{2K_0}}\right\rfloor
 \le \sum_{n\le N}f(n)
\]
for large \(N\), and therefore
\begin{equation}\label{eq:M-bound}
 |Y_N|\ll \frac{1}{\sqrt{N}}\sum_{n\le N}f(n)\ll \sqrt{N},
\end{equation}
Here the last estimate of (\ref{eq:M-bound}) uses $\sum_{n\le N}f(n)\le 2N$ since otherwise our theorem holds trivially.

Now consider two distinct left vertices \(n,n'\in X_N\).  A common neighbour \(w\) gives
\[
 n-w=m^2,
 \qquad
 n'-w={m'}^2
\]
for positive integers \(m\ne m'\).  Suppose, without loss of generality, that \(m>m'\).  Then
\begin{equation}\label{eq:difference-squares}
 |n-n'|=m^2-{m'}^2=(m-m')(m+m').
\end{equation}
Every common neighbour therefore gives a positive divisor
\(
 a=m-m'
\)
of \(h=|n-n'|\).  Once \(a\) is fixed, the other factor is \(b=h/a\), and then
\[
 m=\frac{a+b}{2},
 \qquad
 m'=\frac{b-a}{2}.
\]
Thus \(a\) determines the pair \((m,m')\), if the displayed quantities are positive integers, and then \(w=n-m^2\) is determined as well.  The parity and positivity conditions can only reduce the number of possibilities.  Consequently,
\begin{equation}\label{eq:codegree}
 |\Gamma(n)\cap\Gamma(n')|
 \le\tau(|n-n'|)
 \le T(N),
\end{equation}
because \(1\le|n-n'|\le N\).

Define the first excess moment of the graph by
\[
 a_N=\sum_{n\in X_N}\big(|\Gamma(n)|-1\big)_+.
\]
For large \(N\), any \(n\in X_N\) satisfies \(f(n)\ge1\).  Since \(|\Gamma(n)|=f^*(n)\le f(n)\), we get
\[
 \big(|\Gamma(n)|-1\big)_+\le f(n)-1.
\]
Summing over \(X_N\), we obtain
\begin{equation}\label{eq:a-Delta}
 a_N\le \sum_{n\in X_N}\big(f(n)-1\big)\le \sum_{n\le N}\big(f(n)-1\big)+O(1).
\end{equation}

Apply Lemma~\ref{lem:graph} to \(G_N\), using
\eqref{eq:graph-B}, \eqref{eq:M-bound}, \eqref{eq:codegree}, and \eqref{eq:a-Delta}.  We obtain
\begin{align}\label{eq:almost-final}
 \sum_{N/K_0<n\le N}\binom{f^*(n)}2&=\sum_{n\in X_N}\binom{|\Gamma(n)|}2\nonumber\\
 &\le a_N\sqrt{|Y_N|T(N)}+\frac23|Y_N|\max_{n\in X_N}|\Gamma(n)|\nonumber\\
 &\le \sqrt{T(N)}N^{1/4}\sum_{n\le N}\big(f(n)-1\big)+\frac{N}{(\log N)^{\delta}(\log\log N)^{3/2}}.
\end{align}
Combining \eqref{eq:almost-final} with \eqref{second-moment}, we conclude that
\begin{align*}
\sqrt{T(N)}N^{1/4}\sum_{n\le N}\big(f(n)-1\big)\ge \frac{N\log K_0}{128K_0}-o(N),
\end{align*}
which implies
$$
\sum_{n\le N}f(n)-N\gg \frac{N^{3/4}}{\sqrt{T(N)}},
$$
completing the proof of Theorem \ref{thm:main}.
\end{proof}

\begin{proof}[Proof of Corollary \ref{cor:wigert}]
For every fixed \(\varepsilon>0\),
\begin{equation*}
 \tau(h)
 <2^{(1+\varepsilon)\log h/\log\log h}
\end{equation*}
for every sufficiently large \(h\), see, e.g., Hardy and Wright \cite[Chapter XVIII, \S18.1, Theorem 317, p.~262]{HW}. Therefore,
\begin{equation*}
 T(N)
 \le
 \exp\!\left(
   \big(\log 2+o(1)\big)\frac{\log N}{\log\log N}
 \right),
\end{equation*}
from which our corollary follows.
\end{proof}

\section{Final remarks}
We now discuss the relationship between Problem \ref{pro:2} and Problem \ref{pro:1}. We will see that they are equivalent.

\begin{theorem}\label{thm:3}
If $\mathcal{W}=\{w_n\}_{n=1}^{\infty}$ is an additive complement of the squares that is exact on average, then 
$$
w_n \sim \frac{\pi^2 n^2}{16},  \quad \text{or~equivalently} \quad \mathcal{W}(N)\sim \frac{4}{\pi}\sqrt{N}.
$$
\end{theorem}
\begin{proof}
Throughout the proof, let $N$ be sufficiently large. By (\ref{eq-1216-1}) we have
\begin{align*}
\limsup_{N\rightarrow\infty}\mathcal{W}\big(N\big)/\sqrt{N}\ge \liminf_{N\rightarrow\infty}\mathcal{W}\big(N\big)/\sqrt{N}\ge \frac{4}{\pi}.
\end{align*}
It remains to rule out the following two cases.

{\it Case I.} 
$$
\limsup_{N\rightarrow\infty}\frac{\mathcal{W}(N)}{\sqrt{N}}>\sqrt{6}.
$$
In this case, there are infinitely many positive integers $N$ such that
$
\mathcal{W}(N)>\sqrt{6N}.
$
For these $N$, we clearly have
$$
\sum_{n\le 2N}f(n)\ge \sum_{w\le N,\ m^2\le N}1\ge \sqrt{6N}\cdot \big(1+o(1)\big)\sqrt{N}>\sqrt{5}N,
$$
from which it follows immediately that
$$
\sum_{n\le 2N}f(n)-2N\ge (\sqrt{5}-2)N.
$$
This contradicts the assumption that $\mathcal{W}$ is exact on average.

{\it Case II.} 
\begin{align*}
\frac{4}{\pi}<\limsup_{N\rightarrow\infty}\frac{\mathcal{W}(N)}{\sqrt{N}}=\gamma\le \sqrt{6}.
\end{align*}
In this case, set 
$$
\delta=\frac{\gamma+4/\pi}{2}.
$$
Then there are infinitely many positive integers $N$ such that
\begin{align}\label{eq-1215-1}
\mathcal{W}(N)>\delta \sqrt{N}.
\end{align} 
For these $N$, let 
$$
\delta_1=\frac{\delta+4/\pi}{2} \quad \text{and} \quad \varepsilon=\left(\frac{\delta-\delta_1}{4}\right)^2.
$$ 
We further split {\it Case II} into two subcases.

{\it Subcase II-1.}  
\begin{align}\label{eq-1215-2}
\mathcal{W}\Big((1-\varepsilon)N\Big)\le\delta_1 \sqrt{N}.
\end{align} 
In this subcase we have 
\begin{align}\label{eq-1220-1}
\mathcal{W}(N)-\mathcal{W}\Big((1-\varepsilon)N\Big)>\big(\delta-\delta_1\big)\sqrt{N}
\end{align} 
by (\ref{eq-1215-1}) and (\ref{eq-1215-2}). Hence, from (\ref{eq-1220-1}) we get
\begin{align*}
\sum_{(1-\varepsilon)N<n\le (1+\varepsilon)N}f(n)&\ge \sum_{\substack{(1-\varepsilon)N<w\le N\\ m^2\le \varepsilon N}}1\\
& \ge \big(\delta-\delta_1\big)\sqrt{N}\cdot\Big(1+o(1)\Big)\sqrt{\varepsilon N}\\
&>\frac{(\delta-\delta_1)^2}{5}N.
\end{align*}
Therefore, it follows that
\begin{align*}
\sum_{n\le (1+\varepsilon)N}f(n)-(1+\varepsilon)N&\ge \sum_{(1-\varepsilon)N<n\le (1+\varepsilon)N}f(n)-2\varepsilon N+O(1)\\
&>\frac{(\delta-\delta_1)^2}{5}N- 2\left(\frac{\delta-\delta_1}{4}\right)^2N+O(1)\\
&>\frac{(\delta-\delta_1)^2}{20}N.
\end{align*}
Hence, in {\it Subcase II-1}, $\mathcal{W}$ cannot be exact on average. 

{\it Subcase II-2.}  
\begin{align}\label{eq-1215-3}
\mathcal{W}\Big((1-\varepsilon)N\Big)>\delta_1 \sqrt{N}.
\end{align} 
In this subcase, we first note that
\begin{align*}
\sum_{n\le N}f(n)=\sum_{w+m^2\le N}1=\sum_{m<\sqrt{N}}\sum_{w\le N-m^2}1=\sum_{m<\sqrt{N}}\mathcal{W}\big(N-m^2\big).
\end{align*} 
We now split the sum according as $m\le \sqrt{\varepsilon N}$ or $m>\sqrt{\varepsilon N}$. Then
\begin{align}\label{eq-1215-4}
\sum_{n\le N}f(n)=\sum_{m\le \sqrt{\varepsilon N}}\mathcal{W}\big(N-m^2\big)+\sum_{\sqrt{\varepsilon N}< m<\sqrt{N}}\mathcal{W}\big(N-m^2\big).
\end{align} 
For $m\le \sqrt{\varepsilon N}$ we have $N-m^2\ge (1-\varepsilon)N$, and hence by (\ref{eq-1215-3}) we get
\begin{align}\label{eq-1215-5}
\mathcal{W}\big(N-m^2\big)>\delta_1 \sqrt{N}.
\end{align}
By (\ref{eq-1216-1}), (\ref{eq-1215-4}) and (\ref{eq-1215-5}) we obtain 
\begin{align}\label{eq-1215-6}
\sum_{n\le N}f(n)&>\delta_1 \sqrt{N}\sum_{m\le \sqrt{\varepsilon N}}1+\sum_{\sqrt{\varepsilon N}< m<\sqrt{N}}\mathcal{W}\big(N-m^2\big)\nonumber\\
&\ge \big(\delta_1-4/\pi\big) \sqrt{\varepsilon}N +\sum_{m<\sqrt{N}}\Big(4/\pi+o(1)\Big)\sqrt{N-m^2}+O\big(\sqrt{N}\big).
\end{align}
A routine partial summation gives
\begin{align}\label{eq-1215-7}
\sum_{m<\sqrt{N}}\Big(4/\pi+o(1)\Big)\sqrt{N-m^2}\sim N \quad \text{as }N\rightarrow\infty.
\end{align}
Inserting (\ref{eq-1215-7}) into (\ref{eq-1215-6}), we then obtain
\begin{align*}
\sum_{n\le N}f(n)-N&\ge \Big(\delta_1-4/\pi+o(1)\Big) \sqrt{\varepsilon}N\\
&> \frac{\big(\delta-\delta_1\big)\big(\delta_1-4/\pi\big)}{5}N.
\end{align*} 
Hence, in {\it Subcase II-2}, $\mathcal{W}$ cannot be exact on average.
Thus,
\[
 \limsup_{N\to\infty}\frac{\mathcal{W}(N)}{\sqrt N}\le \frac{4}{\pi}.
\]
Together with (\ref{eq-1216-1}), this gives $\mathcal{W}(N)\sim \frac{4}{\pi}\sqrt N$, which is equivalent to $w_n\sim \frac{\pi^2}{16}n^2$.
\end{proof}

\section*{Acknowledgments}
We thank Yong-Gao Chen and Imre Z. Ruzsa for their helpful conversations.

C.S. was supported by NKFIH Grants Nos. K129335, K146387, and KKP 144059.


\begin{thebibliography}{1}
\bibitem{Abbott} 
H. L. Abbott, {\it On the additive completion of sets of integers,} J. Number Theory {\bf 17} (1983), 135--143.

\bibitem{Balasubramanian} 
R. Balasubramanian, {\it On the additive completion of squares,} J. Number Theory {\bf29} (1988), 10--12.

\bibitem{B-Rama} 
R. Balasubramanian and D. S. Ramana, {\it Additive complements of the squares,} C. R. Math. Acad. Sci. Soc. R. Can. {\bf23} (2001), 6--11.

\bibitem{Balasubramanian-2}
R. Balasubramanian and K. Soundararajan, {\it On the additive completion of squares, II,} J. Number Theory {\bf40} (1992), 127--129.

\bibitem{Bloom}
T. F. Bloom, {\it Erd\H os Problem $\#$33,} https://www.erdosproblems.com/33, accessed 2025-12-20.

\bibitem{Bloom2}
T. F. Bloom, {\it Erd\H os Problem $\#$221,} https://www.erdosproblems.com/221, accessed 2025-12-20.

\bibitem{Chang}
M.-C. Chang, {\it Product sets of arithmetic progressions,} Unpublished manuscript. 

\bibitem{Chen-Fang} 
Y.-G. Chen and J.-H. Fang, Additive complements of the squares,  J. Number Theory {\bf 180} (2017), 410--422.

\bibitem{Cilleruelo} 
J. Cilleruelo, {\it The additive completion of $k$-th powers,} J. Number Theory {\bf44} (1993), 237--243.

\bibitem{Ding}
Y. Ding, {\it Green’s problem on additive complements of the squares,} C. R. Math. Acad. Sci. Paris {\bf358} (2020), 897--900.

\bibitem{DSWX}
Y. Ding, Y.-C. Sun, L.-Y. Wang and Y. Xia, {\it A note on additive complements of the squares,} Discrete Math. {\bf349} (2026),  Paper 114763, 8 pp.

\bibitem{Donagi} 
R. Donagi and M. Herzog, {\it On the additive completion of polynomial sets of integers,} J. Number Theory {\bf3} (1971), 150--154.

\bibitem{Elekes-Ruzsa} 
Gy. Elekes and I. Z. Ruzsa, {\it Few sums, many products,} Studia Sci. Math. Hungar. {\bf40} (2003), no. 3,
301--308. 

\bibitem{Erdos-New}
P. Erd\H os, {\it An asymptotic inequality in the theory of numbers,} Vestnik Leningrad. Univ. {\bf15} (1960), 41--49.


\bibitem{Erdos} 
P. Erd\H os, {\it Problems and results in additive number theory,} in: Colloque sur la Th\'eorie des Nombres, Bruxelles, 1955, George Thone, Li\`ege Masson and Cie, Paris, 1956, pp. 127--137.



\bibitem{Erdos2} 
P. Erd\H os, {\it Problem 33,} Proc. Number Theory Conf., Boulder, Colorado, 1963.

\bibitem{Ford}
K. Ford, {\it The distribution of integers with a divisor in a given interval,} Ann. of Math. (2) {\bf 168} (2008), 367--433. 

\bibitem{Habsieger} 
L. Habsieger, {\it On the additive completion of polynomial sets,} J. Number Theory {\bf51} (1995), 130--135.

\bibitem{HW}
G.~H. Hardy and E.~M. Wright,
\emph{An Introduction to the Theory of Numbers},
4th ed., Oxford University Press, 1960 (corrected reprints).

\bibitem{Moser} 
L. Moser, {\it On the Additive Completion of Sets of Integers,} Proceedings of Symposia in Pure Mathematics, vol. VIII, Amer. Math. Soc., Providence, RI, 1965, pp. 175--180.

\bibitem{Ramana} 
D. S. Ramana, {\it Some Topics in Analytic Number Theory,} PhD thesis, University of Madras, May 2000.

\bibitem{Ramana-2}
D. S. Ramana, {\it A report on additive complements of the squares,} in: Number Theory and Discrete Mathematics, Chandigarh, 2000, in: Trends Math., Birkh\"auser, Basel, 2002, pp. 161--167.

\bibitem{Ruzsa2}
I. Z. Ruzsa, {\it Additive completion of lacunary sequences,} Combinatorica, {\bf21} (2001), 279--291.

\bibitem{Ruzsa}
I. Z. Ruzsa, {\it On a problem of P. Erd\H os,} Canad. Math. Bull. {\bf15} (1972), 309--310.

\bibitem{Tenenbaum}
G. Tenenbaum, {\it Introduction to analytic and probabilistic number theory, third ed.,} Graduate Studies in
Mathematics, vol. 163, American Mathematical Society, Providence, RI, 2015, Translated from the
2008 French edition by Patrick D. F. Ion.

\bibitem{Xu-Zhou}
Max W. Xu and Y. Zhou, {\it On product sets of arithmetic progressions,} Discrete Anal. (2023), Paper No. 10, 31 pp.  

\end{thebibliography}
\end{document}